\documentclass[reqno,b5paper]{amsart}

\usepackage{amsmath}
\usepackage{amssymb}
\usepackage{amsthm}
\usepackage{enumerate}
\usepackage[mathscr]{eucal}
\usepackage{eqlist}

\usepackage[english]{babel}

\theoremstyle{plain}
\newtheorem{theorem}{Theorem}[section]
\newtheorem{prop}[theorem]{Proposition}
\newtheorem{lemma}{Lemma}[section]

\theoremstyle{definition}
\newtheorem{definition}{Definition}[section]

\theoremstyle{remark}

\numberwithin{equation}{section}

\begin{document}

\title[Large and moderate deviation for averaged stochastic algorithm]{Large and moderate deviation principles for averaged stochastic approximation method for the estimation of a regression function}

\author{Yousri Slaoui*}
\newcommand{\acr}{\newline\indent}
\address{\llap{*\,}Universit\'e de Poitiers \acr
Laboratoire de Math\'ematiques et Applications \acr
11 Boulevard Marie et Pierre Curie\acr
86962 Futuroscope Chasseneuil\acr
France}
\email{Yousri.Slaoui@math.univ-poitiers.fr}



\maketitle

\begin{abstract}
In this paper we prove large deviations principles
for the averaged stochastic approximation method for the estimation of a regression function introduced by A. Mokkadem et al. [Revisiting R\'ev\'esz's stochastic approximation method for the estimation of a regression function, ALEA Lat. Amm. J. Probab. Math. Stat. {\bf 6} (2009), 63--114]. We show that the averaged stochastic approximation algorithm constructed using the weight sequence which minimize the asymptotic variance gives the same pointwise LDP as the Nadaraya-Watson kernel estimator. Moreover, we give a moderate deviations principle for these estimators. It turns
out that the rate function obtained in the moderate deviations
principle for the averaged stochastic approximation algorithm constructed using the weight sequence which minimize the asymptotic variance is larger
than the one obtained for the Nadaraya-Watson estimator and the one obtained for the semi-recursive estimator.
\end{abstract}

\section{Introduction}
Let $\left(X,Y\right), \left(X_1,Y_1\right),\ldots ,\left(X_n,Y_n\right)$ be independent, identically distributed pairs of random variables with joint density function $g\left(x,y\right)$, and let $f$ denote the probability density of $X$. In order to construct a stochastic
algorithm for the estimation of the regression function $r:x\mapsto
\mathbb{E}\left(Y|X=x\right)$ at a point $x$ such that $f(x)\neq 0$, A. Mokkadem et al. \cite{Mok09} defines an algorithm, which approximates the zero of the function
$h~:y\mapsto f(x)r(x)-f(x)y$. Following Robbins-Monro's procedure,
this algorithm is defined by setting $r_0(x)\in\mathbb{R}$ and, for $n\geq
1$,
\begin{eqnarray*}
r_n(x)  =  r_{n-1}(x)+\gamma_n{\mathcal W}_n(x)
\end{eqnarray*}
where ${\mathcal W}_n(x)$ is an ``observation'' of the function $h$ at the
point $r_{n-1}(x)$. To define ${\mathcal W}_n(x)$, A. Mokkadem et al \cite{Mok09} follow the approach of P. R\'ev\'esz (\cite{Rev73},~\cite{Rev77}) and A. B. Tsybakov~\cite{Tsy90}, and introduces a kernel $K$ (that is, a function satisfying $\int_\mathbb{R}
K(x)dx=1$) and a bandwidth $\left(h_n\right)$ (that is, a sequence of positive
real numbers that goes to zero), and sets
\begin{eqnarray*}
\mathcal {W}_n(x)=h_n^{-1}Y_nK(h_n^{-1}[x-X_n])-h_n^{-1}K(h_n^{-1}[x-X_n])r_{n-1}(x).
\end{eqnarray*}
Then, the estimator $r_n$ can be rewritten as
\begin{eqnarray}\label{algorevgen}
r_n(x)  =  
\left(1-\gamma_nh_n^{-1}K\left(\frac{x-X_n}{h_n}\right)\right)r_{n-1}(x)
+\gamma_nh_n^{-1}Y_nK\left(\frac{x-X_n}{h_n}\right).
\end{eqnarray}
Now, let the stepsize in (\ref{algorevgen}) satisfy
$\lim_{n\to \infty}n\gamma_n=\infty$, and let $(q_n)$ be a positive sequence
such that $\sum q_n=\infty$. The averaged stochastic approximation algorithm for the estimation of a regression function is defined by setting
\begin{equation}
\label{algo:avergrev}
\overline r_n(x) = \frac{1}{\sum_{k=1}^nq_k}\sum_{k=1}^nq_kr_k(x)
\end{equation}
(where the $r_k(x)$ are given by the algorithm (\ref{algorevgen})).

Recently, large and moderate deviations results have been proved for the
well-known nonrecursive Nadaraya-Watson's kernel regression estimator, first by Louani (1999), and then by C. Joutard \cite{Jot06}. A. Mokkadem et al \cite{Mok08} show that the rate function obtained in the moderate deviations principle for the semi-recursive estimator is larger than the one obtained for the Nadaraya-Watson estimator.\\

Let us first recall that a $\mathbb{R}^m$-valued sequence $\left(Z_n\right)_{n\geq 1}$ satisfies a large deviations principle (LDP) with speed $\left(\nu_n\right)$ and good rate function $I$ if :

\begin{enumerate}
\item $\left(\nu_n\right)$ is a positive sequence such that $\lim_{n\to \infty}\nu_n=\infty$;
\item $I:\mathbb{R}^m\to \left[0,\infty\right]$ has compact level sets;
\item for every borel set $B\subset \mathbb{R}^m$,
\begin{eqnarray*}
-\inf_{x\in \overset{\circ}{B}}I\left(x\right)&\leq & \liminf_{n\to \infty} \nu_n^{-1}\log \mathbb{P}\left[Z_n\in B\right]\\
&\leq & \limsup_{n\to \infty} \nu_n^{-1}\log \mathbb{P}\left[Z_n\in B\right]\leq -\inf_{x\in \overline{B}}I\left(x\right),
\end{eqnarray*}
 where $\overset{\circ}{B}$ and $\overline{B}$ denote the interior and the closure of $B$ respectively. Moreover, let $\left(v_n\right)$ be a nonrandom sequence that goes to infinity; if $\left(v_nZ_n\right)$ satisfies a LDP,
then $\left(Z_n\right)$ is said to satisfy a moderate deviations principle (MDP).
\end{enumerate}

The first aim of this paper is to establish pointwise LDP for the averaged stochastic approximation algorithm~(\ref{algo:avergrev}). It turns out that the rate function depend on the bandwidths $\left(h_n\right)$ and on the weight $\left(q_n\right)$.\\

We show that using the bandwidths $\left(h_n\right)\equiv \left(cn^{-a}\right)$ with $c>0$ and $a\in\left]1-\alpha,\left(4\alpha-3\right)/2\right[$ (with $\alpha\in ]\frac{3}{4},1]$), and the weight $\left(q_n\right)=\left(c^{\prime}n^{-q}\right)$ with $c^{\prime}>0$ and $q<\min\left\{1-2a,\left(1+a\right)/2\right\}$, the sequence $\left(\overline r_{n}\left(x\right)-r\left(x\right)\right)$ satisfies a LDP with speed $\left(nh_n\right)$ and the rate function defined as follows:

\begin{eqnarray*}
I_{a,q,x}\left(t\right)=\sup_{u\in \mathbb{R}}\left\{ut-\psi_{a,q,x}\left(u\right)\right\},
\end{eqnarray*}
which is the Fenchel-Legendre transform of the function $\psi_{a,q,x}$ defined as follows:
\begin{eqnarray}\label{eq:psiaqx}
\psi_{a,q,x}\left(u\right)=\left(1-q\right)\int_{\left[0,1\right]\times\mathbb{R}^2}s^{-a}\left(e^{us^{a-q}K\left(z\right)\frac{\left(y-r\left(x\right)\right)}{f\left(x\right)}}-1\right)g\left(x,y\right)dsdzdy.
\end{eqnarray}
Noting that, in the special case $\left(q_n\right)=\left(h_n\right)$, which is the case when the weight $\left(q_n\right)$ minimizes the asymptotic variance of $\overline r_{n}$ (see A. Mokkadem et al., \cite{Mok09}), we obtain the same rate function for the pointwise LDP as the one obtained for the Nadaraya-Watson estimator (see D. Louani, \cite{Lou99}). \\

Our second aim is to provide pointwise MDP for the averaged stochastic approximation algorithm~(\ref{algo:avergrev}). In this case, we consider more general weight sequence defined as $q_n=\gamma\left(n\right)$ for all $n$, where $\gamma$ is a regularly function with exponent $\left(-q\right)$, $q<\min\left\{1-2a,\left(1+a\right)/2\right\}$.\\

For any positive sequence $\left(v_n\right)$ satisfying 
\begin{eqnarray}\label{cond:mdp}
\lim_{n\to \infty}v_n=\infty, \quad \quad \lim_{n\to \infty}\frac{v_n^2}{nh_n}=0 \quad\mbox{and} \quad \lim_{n\to \infty}v_nh_n^2=0
\end{eqnarray}
and general bandwidths $\left(h_n\right)$, we prove that the sequence
\begin{eqnarray*}
v_n\left(\overline r_n\left(x\right)-r\left(x\right)\right)
\end{eqnarray*} 
satisfies a LDP of speed $\left(nh_n/v_n^2\right)$ and good rate function $J_{a,q,x}:\mathbb{R}\to \mathbb{R}$ defined by 
\begin{eqnarray}\label{MDPalgoreg}
J_{a,q,x}\left(t\right)=\frac{1+a-2q}{\left(1-q\right)^2}\frac{f\left(x\right)}{Var\left[Y\vert X=x\right]\int_{\mathbb{R}}K^2\left(z\right)dz}\frac{t^2}{2}.
\end{eqnarray}
Let us point out that when the weight $\left(q_n\right)$ is chosen to be a regularly varying function with exponent $\left(-a\right)$ (e.g. $\left(q_n\right)=\left(h_n\right)$), which is the case when the weight $\left(q_n\right)$ minimizes the asymptotic variance of $\overline r_{n}$ (see A. Mokkadem et al., \cite{Mok09}), the factor $\left(1+a-2q\right)/\left(1-q\right)^2$ which is present in~(\ref{MDPalgoreg}) can be reduced to $1/(1-a)$, and then we can write
\begin{eqnarray}\label{MDPalgoreg1}
J_{a,x}\left(t\right)=\frac{1}{\left(1-a\right)}\frac{f\left(x\right)}{Var\left[Y\vert X=x\right]\int_{\mathbb{R}}K^2\left(z\right)dz}\frac{t^2}{2}.
\end{eqnarray}

Moreover, D. Louani \cite{Lou99} establish the moderate deviations behaviour for the Nadaraya-Watson (\cite{Nad64},~\cite{Wat64}) estimator defined as
\begin{eqnarray}\label{eq:NW}
\widehat{r}_n\left(x\right)=\left\{\begin{array}{ll}
\frac{\widehat{m}_n\left(x\right)}{\widehat{f}_n\left(x\right)} & \quad if  \quad \widehat{f}_n\left(x\right)\not=0\\
0 \quad &\quad  otherwise,
\end{array}
\right.
\end{eqnarray}
where
\begin{eqnarray*}
\widehat{m}_n\left(x\right)=\frac{1}{nh_n}\sum_{i=1}^nY_iK\left(\frac{x-X_i}{h_n}\right)\quad \mbox{and}\quad \widehat{f}_n\left(x\right)=\frac{1}{nh_n}\sum_{i=1}^nK\left(\frac{x-X_i}{h_n}\right).
\end{eqnarray*}

They prove that, for any positive sequence $\left(v_n\right)$ satisfying~(\ref{cond:mdp}), the sequence $v_n\left(\widehat{r}_n\left(x\right)-r\left(x\right)\right)$ satisfies a LDP with speed $\left(nh_n/v_n^2\right)$ and good rate function 
$\widehat{J}_{x}:\mathbb{R}\to \mathbb{R}$ defined by 
\begin{eqnarray}\label{MDPalgoregNW}
\widehat{J}_{x}\left(t\right)=\frac{f\left(x\right)}{Var\left[Y\vert X=x\right]\int_{\mathbb{R}}K^2\left(z\right)dz}\frac{t^2}{2}.
\end{eqnarray}

Recently, A. Mokkadem et al \cite{Mok08} establish the moderate deviations behaviour for the semi-recursive version of the Nadaraya-Watson estimator defined as
\begin{eqnarray}\label{eq:semi_rec}
\tilde{r}_n\left(x\right)=\left\{\begin{array}{ll}
\frac{\tilde{m}_n\left(x\right)}{\tilde{f}_n\left(x\right)} & \quad if  \quad \tilde{f}_n\left(x\right)\not=0\\
0 \quad &\quad  otherwise,
\end{array}
\right.
\end{eqnarray}
where
\begin{eqnarray*}
\tilde{m}_n\left(x\right)=\frac{1}{n}\sum_{i=1}^n\frac{Y_i}{h_i}K\left(\frac{x-X_i}{h_i}\right)\quad \mbox{and}\quad \tilde{f}_n\left(x\right)=\frac{1}{n}\sum_{i=1}^n\frac{1}{h_i}K\left(\frac{x-X_i}{h_i}\right).
\end{eqnarray*}

They prove that, for any positive sequence $\left(v_n\right)$ satisfying~(\ref{cond:mdp}), the sequence $v_n\left(\tilde{r}_n\left(x\right)-r\left(x\right)\right)$ satisfies a LDP with speed $\left(nh_n/v_n^2\right)$ and good rate function 
$\tilde{J}_{a,x}:\mathbb{R}\to \mathbb{R}$ defined by 
\begin{eqnarray}\label{MDPalgoregsemirec}
\tilde{J}_{a,x}\left(t\right)=\left(1+a\right)\frac{f\left(x\right)}{Var\left[Y\vert X=x\right]\int_{\mathbb{R}}K^2\left(z\right)dz}\frac{t^2}{2}.
\end{eqnarray}

Then, it follows from~(\ref{MDPalgoreg1}),~(\ref{MDPalgoregNW}) and~(\ref{MDPalgoregsemirec}), that the rate function obtained in the MDP of $\overline r_n$ defined with a weight $\left(q_n\right)$ minimizing the asymptotic variance of $\overline r_n$ (e.g. $\left(q_n\right)=\left(h_n\right)$) is larger than the one obtained for the Nadaraya-Watson kernel estimator~(\ref{eq:NW}) and than the one obtained for the semi-recursive kernel estimator~(\ref{eq:semi_rec}); this means that
the averaged stochastic approximation algorithm $\overline r_n(x)$ defined with a weight $\left(q_n\right)$, which is chosen to be a regularly varying function with exponent $\left(-a\right)$ (e.g. $\left(q_n\right)=\left(h_n\right)$) is more concentrated around $r(x)$ than the two others estimators (Nadaraya-Watson~(\ref{eq:NW}) and semi-recursive~(\ref{eq:semi_rec})).

\section{Assumptions and main results} \label{asection 2}
Let us first define the class of positive sequences that will be used
in the statement of our assumptions.

\begin{definition}
Let $\gamma \in \mathbb{R} $ and $\left(v_n\right)_{n\geq 1}$ be a nonrandom
positive sequence. We say that $\left(v_n\right) \in \mathcal{GS}\left(\gamma
\right)$ if
\begin{eqnarray}\label{aGS}
\lim_{n \to \infty} n\left[1-\frac{v_{n-1}}{v_{n}}\right]=\gamma.
\end{eqnarray}
\end{definition}
Condition~\eqref{aGS} was introduced by J. Galambos and E. Seneta \cite{Ga73} to define
regularly varying sequences (see also R. Bojanic and E. Seneta \cite{Bo73}); it was used in 
A. Mokkadem and M. Pelletier \cite{Mok07} in the context of stochastic approximation algorithms.
Typical sequences in $\mathcal{GS}\left(\gamma\right)$ are, for $b\in \mathbb{R}$,
$n^{\gamma}\left(\log n\right)^{b}$, $n^{\gamma}\left(\log \log n\right)^{b}$, and
so on.\\

Let $g\left(s,t\right)$ denote the density of the couple
$\left(X,Y\right)$ (in particular
$f\left(x\right)=\int_{\mathbb{R}}g\left(x,t\right)dt$), and set
$a\left(x\right)=r\left(x\right)f\left(x\right)$.

\subsection{Pointwise LDP for the averaged stochastic approximation algorithm~(\ref{algo:avergrev})}
To establish pointwise LDP for $\overline r_n$, we need the following assumptions.
\begin{itemize}
\item[(L1)] $K:\mathbb{R}\rightarrow \mathbb{R}$ is a nonnegative, continuous, bounded function satisfying $\int_{\mathbb{R}}K\left( z\right) dz=1$, $\int_{\mathbb{R}}zK\left( z\right) dz=0$ and $\int_{\mathbb{R}}z^2K\left( z\right) dz<\infty$. 
\item[(L2)] i) $\left(\gamma_{n}\right)=\mathcal{GS}\left(-\alpha\right)$ with $\alpha \in ]\frac{3}{4},1]$; $\lim_{n\to \infty}n\gamma_n\left(\ln \left(\sum_{k=1}^n
\gamma_k\right)\right)^{-1}=\infty$. \\
  $ii)$ $\left(h_n\right)=\left(cn^{-a}\right)$ with $a\in\left]1-\alpha,\left(4\alpha-3\right)/2\right[$ and $c>0$.\\
 $iii)$ $\left(q_n\right) = \left(c^{\prime}n^{-q}\right)$ with $q<\min\left\{1-2a,\left(1+a\right)/2\right\}$ and $c^{\prime}>0$.\\
\item[(L3)] $i)$ $g\left(s,t\right)$ is two times continuously differentiable with
respect to $s$.\\
 $ii)$ For $q\in\left\{0,1,2\right\}$, $s \mapsto
\int_{\mathbb{R}}t^qg\left(s,t\right)dt$ is a bounded function continuous at
$s=x$.\\ 
 For $q\in \left[2,3\right]$, $s \mapsto
\int_{\mathbb{R}}\left|t\right|^qg\left(s,t\right)dt$ is a bounded function.\\
$iii)$ For $q\in \left\{0,1\right\}$,
$\int_{\mathbb{R}}\left|t\right|^q\left|\frac{\partial g}{\partial
x}\left(x,t\right)\right|dt<\infty$, and $s\mapsto
\int_{\mathbb{R}}t^q\frac{\partial^2 g}{\partial s^2}\left(s,t\right)dt$ is a
bounded function continuous at $s=x$.\\
\item[(L4)] For any $u\in \mathbb{R}$, $t\to \int_{\mathbb{R}}\exp\left(uy\right)g\left(t,y\right)dy$ is continuous at $x$ and bounded.
\end{itemize}
The proof of the following comment is given in A. Mokkadem et al. \cite{Mok08}.
\paragraph{Comment}
Notice that $\left(L4\right)$ implies that $\forall m\geq 0, \forall \rho \geq 0$ 
\begin{eqnarray}\label{eq:bounded}
\quad \mbox{the function }\quad t\mapsto \int_{\mathbb{R}}\left|y\right|^m\exp\left(\rho\left|y\right|\right)g\left(t,y\right)dy\quad \mbox{is bounded}.
\end{eqnarray}
Before stating our results, we set
\begin{eqnarray*}
S_{+}=\left\{x\in \mathbb{R}; K\left(x\right)>0\right\}\quad \mbox{and}\quad S_{-}=\left\{x\in \mathbb{R}; K\left(x\right)<0\right\}
\end{eqnarray*}
and for fixed $x\in \mathbb{R}$
\begin{eqnarray*}
T_{+}=\left\{y\in \mathbb{R}; y-r\left(x\right)>0\right\}\quad \mbox{and}\quad T_{-}=\left\{y\in \mathbb{R}; y-r\left(x\right)<0\right\}
\end{eqnarray*}
Moreover, we set
\begin{eqnarray*}
O_+=\left(S_+\cap T_+\right) \cup \left(S_-\cap T_-\right)\quad \mbox{and} \quad O_-=\left(S_+\cap T_-\right) \cup \left(S_-\cap T_+\right)
\end{eqnarray*} 

The following proposition gives the properties of the functions $\psi_{a,q,x}$ and $I_{a,q,x}$; in particular, the behaviour of the rate function $I_{a,q,x}$.
\begin{prop}[Properties of $\psi_{a,q,x}$ and $I_{a,q,x}$]\label{prop:convexe}
$ $\\
Let $\lambda$ be the Lebesgue measure on $\mathbb{R}$ and let Assumptions $\left(L1\right)$ and $\left(L4\right)$ hold.
\begin{enumerate}
\item[(i)] $\psi_{a,q,x}$ is strictly convex, twice continuously differentiable on $\mathbb{R}$, and $I_{a,q,x}$ is a good rate function on $\mathbb{R}$.
\item[(ii)] If $\lambda\left(O_-\right)=0$, $I_{a,q,x}\left(t\right)=+\infty$, when $t<0$, and
\begin{eqnarray*}
I_{a,q,x}\left(0\right)&=&\left\{\begin{array}{lllll}
\left(1-q\right)/\left(1-a\right)\lambda\left(S_+\right)f\left(x\right) & \mbox{if} & \lambda\left(S_+\cap T_+\right)>0\\
\left(1-q\right)/\left(1-a\right)\lambda\left(S_-\right)f\left(x\right) & \mbox{if} & \lambda\left(S_-\cap T_-\right)>0\\
\end{array}
\right.\\
\end{eqnarray*}
$I_{a,q,x}$ is strictly convex on $\mathbb{R}$ and continuous on $]0,+\infty[$, and for any $t>0$
\begin{eqnarray}\label{eq:taurev}
I_{a,q,x}\left(t\right)=t\left(\psi_{a,q,x}^{\prime}\right)^{-1}\left(t\right)-\psi_{a,q,x}\left(\left(\psi_{a,q,x}^{\prime}\right)^{-1}\left(t\right)\right),
\end{eqnarray}

\item[(iii)] If $\lambda\left(O_-\right)>0$, then $I_{a,q,x}$ is finite and strictly convex on $\mathbb{R}$ and~(\ref{eq:taurev}) holds for any $t\in \mathbb{R}$.
\end{enumerate}
\end{prop} 
We can now state the LDP for the averaged stochastic approximation algorithm~(\ref{algo:avergrev}).
\begin{theorem}[Pointwise LDP for the averaged stochastic approximation algorithm~\ref{algo:avergrev}]\label{the:LDavergrev} $ $\\
Let Assumptions $\left(L1\right)-\left(L4\right)$ hold. Then, the sequence $\left(\overline r_n\left(x\right)-r\left(x\right)\right)$ satisfies a LDP with speed $\left(nh_n\right)$ and rate function defined as follows:
\begin{eqnarray*}
I_{a,q,x}\left(t\right)=t\left(\psi_{a,q,x}^{\prime}\right)^{-1}\left(t\right)-\psi_{a,q,x}\left(\left(\psi_{a,q,x}^{\prime}\right)^{-1}\left(t\right)\right),
\end{eqnarray*}
where $\psi_{a,q,x}$ is defined in~(\ref{eq:psiaqx}).
\end{theorem}

\subsection{Pointwise MDP for the averaged stochastic approximation algorithm~(\ref{algo:avergrev})}
Let $\left(v_n\right)$ be a positive sequence; we assume that
\begin{itemize} 
\item[(M1)] $K:\mathbb{R}\rightarrow \mathbb{R}$ is a nonnegative, continuous, bounded function satisfying $\int_{\mathbb{R}}K\left( z\right) dz=1$, $\int_{\mathbb{R}}zK\left( z\right) dz=0$ and $\int_{\mathbb{R}}z^2K\left( z\right) dz<\infty$. 
\item[(M2)] i) $\left(\gamma_{n}\right)=\mathcal{GS}\left(-\alpha\right)$ with $\alpha \in ]\frac{3}{4},1]$; $\lim_{n\to \infty}n\gamma_n\left(\ln \left(\sum_{k=1}^n
\gamma_k\right)\right)^{-1}=\infty$. \\
  $ii)$ $\left(h_n\right)=\mathcal{GS}\left(-a\right)$ with $a\in\left]1-\alpha,\left(4\alpha-3\right)/2\right[$.\\
 $iii)$ $\left(q_n\right) = \mathcal{GS}\left(-q\right)$ with $q<\min\left\{1-2a,\left(1+a\right)/2\right\}$.\\
\item[(M3)] $i)$ $g\left(s,t\right)$ is two times continuously differentiable with
respect to $s$.\\
 $ii)$ For $q\in\left\{0,1,2\right\}$, $s \mapsto
\int_{\mathbb{R}}t^qg\left(s,t\right)dt$ is a bounded function continuous at
$s=x$.\\ 
 For $q\in \left[2,3\right]$, $s \mapsto
\int_{\mathbb{R}}\left|t\right|^qg\left(s,t\right)dt$ is a bounded function.\\
$iii)$ For $q\in \left\{0,1\right\}$,
$\int_{\mathbb{R}}\left|t\right|^q\left|\frac{\partial g}{\partial
x}\left(x,t\right)\right|dt<\infty$, and $s\mapsto
\int_{\mathbb{R}}t^q\frac{\partial^2 g}{\partial s^2}\left(s,t\right)dt$ is a
bounded function continuous at $s=x$.\\
\item[(M4)] For any $u\in \mathbb{R}$, $t\to \int_{\mathbb{R}}\exp\left(uy\right)g\left(t,y\right)dy$ is continuous at $x$ and bounded.
\item[(M5)] $i)$ $\lim_{n\to \infty}v_n=\infty$ and $\lim_{n\to \infty}\frac{v_n^2}{nh_n}=0$.\\
$ii)$ $\lim_{n\to \infty}v_nh_n^2=0$
\end{itemize}
The following Theorem gives the pointwise MDP for the averaged stochastic approximation algorithm~(\ref{algo:avergrev}).
\begin{theorem}[Pointwise MDP for the averaged stochastic approximation algorithm~(\ref{algo:avergrev})]\label{the:MDP} $ $\\
Let Assumptions $\left(M1\right)-\left(M5\right) $ hold. Then, the sequence $\left(v_n\left(\overline r_n\left(x\right)-r\left(x\right)\right)\right)$ satisfies a MDP with speed $\left(nh_n/v_n^2\right)$ and  good rate function $J_{a,q,x}$ defined in (\ref{MDPalgoreg}).
\end{theorem}

\section{Proofs}
From now on, we set $n_0\geq 3$ such that $\forall k\geq n_0$, $\gamma_k \leq
\left(2\|f\|_{\infty}\right)^{-1}$ and $\gamma_kh_k^{-1}\|K\|_{\infty}\leq 1$.
Moreover, we introduce 
the following notations:
\begin{eqnarray}
Z_n\left(x\right)&=&h_n^{-1}K\left(\frac{x-X_n}{h_n}\right),\nonumber\\
W_n\left(x\right)&=&h_n^{-1}Y_nK\left(\frac{x-X_n}{h_n}\right),\nonumber\\
\eta_n\left(x\right)&=&\left(Y_n-r\left(x\right)\right)K\left(\frac{x-X_n}{h_n}\right),\label{eq:eta}
\end{eqnarray}
As explained in the introduction, we note that the stochastic approximation
algorithm~\eqref{algorevgen} can be rewritten as:
\begin{eqnarray*}
r_n(x) &
=&\left(1-\gamma_nZ_n\left(x\right)\right)r_{n-1}(x)+\gamma_nW_n\left(x\right)\\
&=&\left(1-\gamma_nf\left(x\right)\right)r_{n-1}(x)+\gamma_n\left(f\left(x\right)-Z_n\left(x\right)\right)r_{n-1}\left(x\right)+\gamma_nW_n\left(x\right).
\end{eqnarray*}
To establish the asymptotic behaviour of $\left(r_n\right)$ and
$\left(\overline r_n\right)$, we introduce the auxiliary stochastic approximation
algorithm defined by setting $\rho_n\left(x\right)=r\left(x\right)$ for all $n\leq
n_0-2$, $\rho_{n_0-1}\left(x\right)=r_{n_0-1}\left(x\right)$, and, for $n\geq n_0$,
\begin{eqnarray*}
\rho_n(x)&=&\left(1-\gamma_nf\left(x\right)\right)\rho_{n-1}(x)+\gamma_n\left(f\left(x\right)-Z_n\left(x\right)\right)r\left(x\right)+\gamma_nW_n\left(x\right).
\end{eqnarray*}
It follows that, for $n\geq n_0$,
\begin{eqnarray*}
\rho_{n}\left(x\right)-\rho_{n-1}\left(x\right)&=&-\gamma_nf\left(x\right)\left[\rho_{n-1}(x)-r\left(x\right)\right]+\gamma_n\left[W_n\left(x\right)-r\left(x\right)Z_n\left(x\right)\right],\\
&=&-\gamma_nf\left(x\right)\left[\rho_{n-1}(x)-r\left(x\right)\right]+\gamma_nh_n^{-1}\eta_n\left(x\right),
\end{eqnarray*}
and thus
\begin{eqnarray*}
\rho_{n-1}(x)-r\left(x\right)&=&\frac{h_n^{-1}}{f\left(x\right)}\eta_n\left(x\right)-\frac{1}{\gamma_nf\left(x\right)}\left[\rho_{n}(x)-\rho_{n-1}(x)\right].
\end{eqnarray*}
Then, we can write that
\begin{eqnarray}\label{rhon}
\overline \rho_n\left(x\right)-r\left(x\right)&=&\frac{1}{\sum_{k=1}^nq_k}\sum_{k=1}^nq_k\left[\rho_k\left(x\right)-r\left(x\right)\right]\nonumber\\
&=&\frac{1}{f\left(x\right)}T_n\left(x\right)-\frac{1}{f\left(x\right)}R_n^{\left(0\right)}\left(x\right)
\end{eqnarray}
with 
\begin{eqnarray*}
T_n\left(x\right)&=&\frac{1}{\sum_{k=1}^nq_k}\sum_{k=n_0-1}^nq_kh_k^{-1}\eta_k\left(x\right),\\
R_n^{\left(0\right)}\left(x\right)&=&\frac{1}{\sum_{k=n_0-1}^nq_k}\sum_{k=1}^n\frac{q_k}{\gamma_{k+1}}\left[\rho_{k+1}(x)-\rho_{k}(x)\right].
\end{eqnarray*}
Moreover, it was showen in A. Mokkadem et al \cite{Mok09}, that under the assumptions $\left(M1\right)-\left(M3\right)$, we have 
\begin{eqnarray}\label{av:Rn0}
\left|R_n^{\left(0\right)}\left(x\right)\right|=o\left(\sqrt{n^{-1}h_n^{-1}}+h_n^{-2}\right) \quad \mbox{a.s.,}
\end{eqnarray} 
then, it follows from~(\ref{rhon}) and~(\ref{av:Rn0}) that
\begin{eqnarray*}
\overline \rho_n\left(x\right)-\mathbb{E}\left[\overline \rho_n\left(x\right)\right]
&=&\frac{1}{f\left(x\right)}\frac{1}{\sum_{k=1}^nq_k}\sum_{k=n_0-1}^nq_kh_k^{-1}\left(\eta_k\left(x\right)-\mathbb{E}\left[\eta_k\left(x\right)\right]\right).
\end{eqnarray*}
Let $\left(\Psi_n\right)$, $\left(B_n\right)$ and $\left(\overline \Delta_n\right)$ be the sequences defined as
\begin{eqnarray*}
\Psi_n\left(x\right)&=&\frac{1}{f\left(x\right)}\frac{1}{\sum_{k=1}^nq_k}\sum_{k=n_0-1}^nq_kh_k^{-1}\left(\eta_k\left(x\right)-\mathbb{E}\left[\eta_k\left(x\right)\right]\right),\\
B_n\left(x\right)&=&\mathbb{E}\left[\overline \rho_n\left(x\right)\right]-r\left(x\right),\\
\overline \Delta_n\left(x\right)&=&\overline r_n\left(x\right)-\overline \rho_n\left(x\right).
\end{eqnarray*}
We have:
\begin{eqnarray}\label{eq:psiB}
\overline r_n\left(x\right)-r\left(x\right)=\Psi_n\left(x\right)+B_n\left(x\right)+\overline \Delta_n\left(x\right).
\end{eqnarray}
Theorems~\ref{the:LDavergrev} and~\ref{the:MDP} are consequences of~(\ref{eq:psiB}) and the following propositions.

\begin{prop}[{Pointwise LDP and MDP for $\left(\Psi_n\right)$}]\label{prop:LMDP}\quad \quad \quad \quad \quad \quad \quad \quad \quad \quad \quad \quad \quad \quad \quad \quad \quad \quad \quad \quad
\begin{enumerate}
\item Under the assumptions $\left(L1\right)-\left(L4\right)$, the sequence $\overline \rho_n\left(x\right)-\mathbb{E}\left[\overline \rho_n\left(x\right)\right]$ satisfies a LDP with speed $\left(nh_n\right)$ and rate function $I_{a,q,x}$.
\item Under the assumptions $\left(M1\right)-\left(M5\right)$, the sequence $\left(v_n\Psi_n\left(x\right)\right)$ satisfies a LDP with speed $\left(nh_n/v_n^2\right)$ and rate function $J_{a,q,x}$.
\end{enumerate}
\end{prop}
\begin{prop}[{Convergence rate of $\left(B_n\right)$}]\label{prop:Bn}
$ $\\
Let Assumptions $\left(M1\right)-\left(M3\right)$ hold. Then
\begin{eqnarray*}
B_n\left(x\right)=O\left(h_n^2\right).
\end{eqnarray*}
\end{prop}
The proof of the following proposition is given in A. Mokkadem et al. \cite{Mok09}.
\begin{prop}[{Convergence rate of $\left(\overline\Delta_n\right)$}]\label{prop:Del}
$ $\\
Let Assumptions $\left(M1\right)-\left(M3\right)$ hold. Then
\begin{eqnarray*}
\overline\Delta_n\left(x\right)=o\left(h_n^2+\frac{1}{\sqrt{nh_n}}\right).
\end{eqnarray*}
\end{prop}

Set $x\in \mathbb{R}$; since the assumptions of Theorems~\ref{the:LDavergrev} guarantee that $\lim_{n\to \infty}B_n\left(x\right)=0$ and $\lim_{n\to \infty}\overline \Delta_n\left(x\right)=0$ Theorem~\ref{the:LDavergrev} is a straightforward consequence of the application of Proposition~\ref{prop:LMDP}. Moreover, under the assumptions of Theorem~\ref{the:MDP}, we have by application of Propostion~\ref{prop:Bn}, $\lim_{n\to \infty}v_nB_n\left(x\right)=0$ and $\lim_{n\to \infty}v_n\overline \Delta_n\left(x\right)=0$; Theorem~\ref{the:MDP} thus straightfully follows from the application of Part 2 of Proposition~\ref{prop:LMDP}.\\
We now state a preliminary lemma, which will be used in the proof of Proposition~\ref{prop:LMDP}. For any $u\in \mathbb{R}$, set
\begin{eqnarray*}
\Lambda_{n,x}\left(u\right)&=&\frac{v_n^2}{nh_n}\log \mathbb{E}\left[\exp\left(unh_n\Psi_n\left(x\right)\right)\right],\\
\Lambda_{x}^{L}\left(u\right)&=&\psi_{a,q,x}\left(u\right),\\
\Lambda_{x}^M\left(u\right)&=&\frac{u^2}{2}\frac{\left(1-q\right)^2}{1+a-2q}\frac{Var\left[Y\vert X=x\right]}{f\left(x\right)}
\int_{\mathbb{R}}K^2\left(z\right)dz.
\end{eqnarray*}
\begin{lemma}[Pointwise convergence of $\Lambda_{n,x}$]\label{lemma:convLam}
$ $\\
For all $u\in \mathbb{R}$
\begin{eqnarray*}
\lim_{n\to \infty}\Lambda_{n,x}\left(u\right)=\Lambda_{x}\left(u\right)
\end{eqnarray*}
where
\begin{eqnarray*}
\Lambda_x\left(u\right)=\left\{\begin{array}{lllll}
\Lambda_x^{L}\left(u\right) & \mbox{when} & v_n\equiv 1\mbox{,}&\left(L1\right)-\left(L4\right)\,\,\mbox{hold} \\
\Lambda_x^{M}\left(u\right) & \mbox{when} & v_n\to \infty \mbox{,}&\left(M1\right)-\left(M4\right)\,\,\mbox{hold} \\
\end{array}
\right.
\end{eqnarray*}
\end{lemma}
Our proofs are now organized as follows: Lemma~\ref{lemma:convLam} is proved in Section~\ref{proof:lemma1}, Proposition~\ref{prop:LMDP} in Section~\ref{proof:prop2} and Proposition~\ref{prop:Bn} in Section~\ref{proof:prop3}.
\subsection{Proof of Lemma~\ref{lemma:convLam}.}\label{proof:lemma1}
\begin{proof}
$ $\\
Set $u\in \mathbb{R}$, $u_n=u/v_n$ and $a_n=nh_n$. We have:
\begin{eqnarray*}
\Lambda_{n,x}\left(u\right)&=&\frac{v_n^2}{a_n}\log \mathbb{E}\left[\exp\left(u_na_n\Psi_n\left(x\right)\right)\right]\\
&=&\frac{v_n^2}{a_n}\log \mathbb{E}\left[\exp\left(\frac{u_n}{f\left(x\right)}\frac{a_n}{\sum_{k=1}^nq_k}\sum_{k=n_0-1}^n\frac{q_k}{h_k}\left(\eta_k\left(x\right)-\mathbb{E}\left[\eta_k\left(x\right)\right]\right)\right)\right]\\
&=&\frac{v_n^2}{a_n}\sum_{k=n_0-1}^n\log \mathbb{E}\left[\exp\left(u_n\frac{a_n}{\sum_{k=1}^nq_k}\frac{q_k}{h_k}\frac{\eta_k\left(x\right)}{f\left(x\right)}\right)\right]\\
&&-\frac{u}{f\left(x\right)} \frac{v_n}{\sum_{k=1}^nq_k}\sum_{k=n_0-1}^n\frac{q_k}{h_k}\mathbb{E}\left[\eta_k\left(x\right)\right].
\end{eqnarray*}
By Taylor expansion, there exists $c_{k,n}$ between $1$ and $\mathbb{E}\left[\exp\left(u_n\frac{q_k}{h_k}\frac{\eta_k\left(x\right)}{f\left(x\right)}\right)\right]$ such that 
\begin{eqnarray*}
\lefteqn{
\log \mathbb{E}\left[\exp\left(u_n\frac{a_n}{\sum_{k=1}^nq_k}\frac{q_k}{h_k}\frac{\eta_k\left(x\right)}{f\left(x\right)}\right)\right]}\\
&&=\mathbb{E}\left[\exp\left(u_n\frac{a_n}{\sum_{k=1}^nq_k}\frac{q_k}{h_k}\frac{\eta_k\left(x\right)}{f\left(x\right)}\right)-1\right]\\
&&-\frac{1}{2c_{k,n}^2}\left(\mathbb{E}\left[\exp\left(u_n\frac{a_n}{\sum_{k=1}^nq_k}\frac{q_k}{h_k}\frac{\eta_k\left(x\right)}{f\left(x\right)}\right)-1\right]\right)^2
\end{eqnarray*}
and $\Lambda_{n,x}$ can be rewriten as
\begin{eqnarray}\label{eq:lamb}
\Lambda_{n,x}\left(u\right)&=&\frac{v_n^2}{a_n}\sum_{k=n_0-1}^n\mathbb{E}\left[\exp\left(u_n\frac{a_n}{\sum_{k=1}^nq_k}\frac{q_k}{h_k}\frac{\eta_k\left(x\right)}{f\left(x\right)}\right)-1\right]
\nonumber\\
&&-\frac{1}{2}\frac{v_n^2}{a_n}\sum_{k=n_0-1}^n
\frac{1}{c_{k,n}^2}\left(\mathbb{E}\left[\exp\left(u_n\frac{a_n}{\sum_{k=1}^nq_k}\frac{q_k}{h_k}\frac{\eta_k\left(x\right)}{f\left(x\right)}\right)-1\right]\right)^2\nonumber\\
&&-\frac{u}{f\left(x\right)} \frac{v_n}{\sum_{k=1}^nq_k}\sum_{k=n_0-1}^n\frac{q_k}{h_k}\mathbb{E}\left[\eta_k\left(x\right)\right].
\end{eqnarray}
Now, let us recall that, if $\left(b_n\right)\in \mathcal{GS}\left(-b^*\right)$ 
with $b^*<1$, then we have, for any fixed $k_0\geq 1$,
\begin{eqnarray}
\lim_{n \to \infty}\frac{nb_n}{\sum_{k=k_0}^nb_k}&=&1-b^*,\label{limbnsum}
\end{eqnarray}
and
\begin{eqnarray}\label{suprapp}
\sup_{k\leq n}\frac{b_n}{b_k}<\infty.
\end{eqnarray}
Moreover, since $\left(q_kh_k^{-1}\right)\in \mathcal{GS}\left(-\left(q-a\right)\right)$, it follows from~(\ref{limbnsum}) that
\begin{eqnarray*}
\left|u_n\frac{a_n}{\sum_{k=1}^nq_k}\frac{q_k}{h_k}\right|&=&O\left(\frac{u}{v_n}\frac{h_n}{h_k}\frac{q_k}{q_n}\right),
\end{eqnarray*}
and from~(\ref{suprapp}) that
\begin{eqnarray*}
\left|u_n\frac{a_n}{\sum_{k=1}^nq_k}\frac{q_k}{h_k}\right|=\left\{\begin{array}{lllll}
O\left(1\right) & \mbox{when} & v_n\equiv 1\\
o\left(1\right) & \mbox{when} & v_n\to \infty\\
\end{array}
\right.
\end{eqnarray*}
and thus, in the both cases, there exists $c>0$ such that
\begin{eqnarray}\label{cte:c}
\left|u_n\frac{a_n}{\sum_{k=1}^nq_k}\frac{q_k}{h_k}\right|&\leq &c.
\end{eqnarray}
{\bf First case: $v_n\to \infty$.} 
$ $\\
A Taylor's expansion implies the existence of $c^{\prime}_{k,n}$ between $0$ and $u_n\frac{a_n}{\sum_{k=1}^nq_k}\frac{q_k}{h_k}\frac{\eta_k\left(x\right)}{f\left(x\right)}$ such that 
\begin{eqnarray*}
\lefteqn{\mathbb{E}\left[\exp\left(u_n\frac{a_n}{\sum_{k=1}^nq_k}\frac{q_k}{h_k}\frac{\eta_k\left(x\right)}{f\left(x\right)}\right)-1\right]}\\
&=&\frac{u_n}{f\left(x\right)}\frac{a_n}{\sum_{k=1}^nq_k}\frac{q_k}{h_k}\mathbb{E}\left[\eta_k\left(x\right)\right]+\frac{1}{2}\left(\frac{u_n}{f\left(x\right)}\frac{a_n}{\sum_{k=1}^nq_k}\frac{q_k}{h_k}\right)^2\mathbb{E}\left[\eta_k^2\left(x\right)\right]
\\
&&+\frac{1}{6}\left(\frac{u_n}{f\left(x\right)}\frac{a_n}{\sum_{k=1}^nq_k}\frac{q_k}{h_k}\right)^3\mathbb{E}\left[\eta_k^3\left(x\right)e^{c^{\prime}_{k,n}}\right].
\end{eqnarray*}
Therefore,
\begin{eqnarray*}
\Lambda_{n,x}\left(u\right)&=&\frac{u^2}{2\left(f\left(x\right)\right)^2}\frac{a_n}{\left(\sum_{k=1}^nq_k\right)^2}\sum_{k=n_0-1}^n\frac{q_k^2}{h_k^2}\mathbb{E}\left[\eta_k^2\left(x\right)\right]\\
&&+\frac{1}{6}\frac{u^2u_n}{\left(f\left(x\right)\right)^3}\frac{a_n^2}{\left(\sum_{k=1}^nq_k\right)^3}\sum_{k=n_0-1}^n\frac{q_k^3}{h_k^3}\mathbb{E}\left[\eta_k^3\left(x\right)\exp\left(c^{\prime}_{k,n}\right)\right]\nonumber\\
&&-\frac{1}{2}\frac{v_n^2}{a_n}\sum_{k=n_0-1}^n
\frac{1}{c_{k,n}^2}\left(\mathbb{E}\left[\exp\left(u_n\frac{a_n}{\sum_{k=1}^nq_k}\frac{q_k}{h_k}\frac{\eta_k\left(x\right)}{f\left(x\right)}\right)-1\right]\right)^2.
\end{eqnarray*}
Let us note that under the assumption $\left(M3\right)$, we have
\begin{eqnarray*}
\mathbb{E}\left[\eta_k^2\left(x\right)\right]=h_kVar\left[Y\vert X=x\right]f\left(x\right)\int_{\mathbb{R}}K^2\left(z\right)dz\left[1+o\left(1\right)\right].
\end{eqnarray*}
Then, it follows that
\begin{eqnarray}\label{eq:Moderate}
\Lambda_{n,x}\left(u\right)&=&\frac{u^2}{2}\frac{a_n}{\left(\sum_{k=1}^nq_k\right)^2}\sum_{k=n_0-1}^n\frac{q_k^2}{h_k}\frac{Var\left[Y\vert X=x\right]}{f\left(x\right)}\int_{\mathbb{R}}K^2\left(z\right)dz\left[1+o\left(1\right)\right]\nonumber\\
&&+R_{n,x}^{\left(1\right)}\left(u\right)-R_{n,x}^{\left(2\right)}\left(u\right),
\end{eqnarray}
with
\begin{eqnarray*}
R_{n,x}^{\left(1\right)}\left(u\right)&=&\frac{1}{6}\frac{u^2u_n}{\left(f\left(x\right)\right)^3}\frac{a_n^2}{\left(\sum_{k=1}^nq_k\right)^3}\sum_{k=n_0-1}^n\frac{q_k^3}{h_k^3}\mathbb{E}\left[\eta_k^3\left(x\right)\exp\left(c^{\prime}_{k,n}\right)\right],\nonumber\\
R_{n,x}^{\left(2\right)}\left(u\right)&=&\frac{1}{2}\frac{v_n^2}{a_n}\sum_{k=n_0-1}^n
\frac{1}{c_{k,n}^2}\left(\mathbb{E}\left[\exp\left(u_n\frac{a_n}{\sum_{k=1}^nq_k}\frac{q_k}{h_k}\frac{\eta_k\left(x\right)}{f\left(x\right)}\right)-1\right]\right)^2.
\end{eqnarray*}
Let us first show that 
\begin{eqnarray*}
\lim_{n\to \infty}\left|R_{n,x}^{\left(1\right)}\left(u\right)\right|=0.
\end{eqnarray*}
In view of $\left(M4\right)$ and~(\ref{cte:c}), we have
\begin{eqnarray}
\lefteqn{\mathbb{E}\left|\eta_k\left(x\right)^3\exp\left(c^{\prime}_{k,n}\right)\right|}\nonumber\\
&\leq&h_k\int_{\mathbb{R}^2}\left|y-r\left(x\right)\right|^{3}K^{3}\left(z\right)\exp\left(\frac{c}{f\left(x\right)}\left|y-r\left(x\right)\right|\left|K\left(z\right)\right|\right)g\left(x-zh_k,y\right)dydz\nonumber\\
&\leq
&4h_k\int_{\mathbb{R}}\exp\left(\frac{c}{f\left(x\right)}\left|r\left(x\right)\right|\left\|K\right\|_{\infty}\right)\left\{\int_{\mathbb{R}}\left|y\right|^{3}\exp\left(\frac{c}{f\left(x\right)}\left|y\right|\left\|K\right\|_{\infty}\right)g\left(x-zh_k,y\right)dy\right.\nonumber\\
&&+\left.\left|r\left(x\right)\right|^{3}\int_{\mathbb{R}}\exp\left(\frac{c}{f\left(x\right)}\left|y\right|\left\|K\right\|_{\infty}\right)g\left(x-zh_k,y\right)dy\right\}K^{3}\left(z\right)dz\nonumber\\
&=&O\left(h_k\right)\label{eta3}.
\end{eqnarray}
Hence, it follows from~(\ref{eta3}) and~(\ref{limbnsum}), that
\begin{eqnarray*}
\lefteqn{\left|\frac{u^2u_n}{\left(f\left(x\right)\right)^2}\frac{a_n^2}{\left(\sum_{k=1}^nq_k\right)^3}\sum_{k=n_0-1}^n\frac{q_k^3}{h_k^3}\mathbb{E}\left[\eta_k^3\left(x\right)e^{c^{\prime}_{k,n}}\right]
\right|}\\
&=& O\left(\frac{1}{v_n}\frac{a_n^2}{\left(\sum_{k=1}^nq_k\right)^3}\sum_{k=n_0-1}^n\frac{q_k^3}{h_k^2}\right)\nonumber\\
&=& O\left(\frac{1}{v_n}\left(\frac{nq_n}{\sum_{k=1}^nq_k}\right)^3\frac{\sum_{k=n_0-1}^nq_k^3h_k^{-2}}{nq_n^3h_n^{-2}}\right)\nonumber\\
&=&O\left(\frac{1}{v_n}\right)
\end{eqnarray*}
which ensures that $\lim_{n\to \infty}\left|R_{n,x}^{\left(1\right)}\left(u\right)\right|=0$.\\
Let us now prove that
\begin{eqnarray*}
\lim_{n\to \infty}\left|R_{n,x}^{\left(2\right)}\left(u\right)\right|=0.
\end{eqnarray*}
Noting that, under the assumption $\left(M3\right)$ we have
\begin{eqnarray*}
\mathbb{E}\left(W_k\left(x\right)\right)&=&a\left(x\right)+\frac{1}{2}h_k^2\int_{\mathbb{R}}y\frac{\partial^2
g}{\partial
x^2}\left(x,y\right)dy\int_{\mathbb{R}}z^2K\left(z\right)dz\left[1+o\left(1\right)\right],\\
\mathbb{E}\left(Z_k\left(x\right)\right)&=&f\left(x\right)+\frac{1}{2}h_k^2\int_{\mathbb{R}}\frac{\partial^2
g}{\partial
x^2}\left(x,y\right)dy\int_{\mathbb{R}}z^2K\left(z\right)dz\left[1+o\left(1\right)\right].
\end{eqnarray*}
Then, it follows from~(\ref{eq:eta}) that
\begin{eqnarray}\label{eq:Eeta}
\mathbb{E}\left[\eta_k\left(x\right)\right]&=&h_k\left[\mathbb{E}\left(W_k\left(x\right)\right)-r\left(x\right)\mathbb{E}\left(Z_k\left(x\right)\right)\right]\nonumber\\
&=&\frac{1}{2}h_k^3\left[\int_{\mathbb{R}}y\frac{\partial^2
g}{\partial
x^2}\left(x,y\right)dy-r\left(x\right)\int_{\mathbb{R}}y\frac{\partial^2
g}{\partial
x^2}\left(x,y\right)dy\right]\nonumber\\
&&\times \int_{\mathbb{R}}z^2K\left(z\right)dz\left[1+o\left(1\right)\right]\nonumber\\
&=&h_k^3m^{\left(2\right)}\left(x\right)f\left(x\right)\left[1+o\left(1\right)\right],
\end{eqnarray}
where,
\begin{eqnarray*}
m^{\left(2\right)}\left(x\right)
=\frac{1}{2f\left(x\right)}\left[\int_{\mathbb{R}}t\frac{\partial^2
g}{\partial x^2}\left(x,t\right)dt-r\left(x\right)\int_{\mathbb{R}}\frac{\partial^2
g}{\partial 
x^2}\left(x,t\right)dt\right]\int_{\mathbb{R}}z^2K\left(z\right)dz.
\end{eqnarray*}
Moreover, in view of~(\ref{limbnsum}) and~(\ref{eq:Eeta}), we have
\begin{eqnarray}\label{eq:R2term2}
\lefteqn{\left|\frac{v_n^2}{a_n}\sum_{k=n_0-1}^n
\frac{1}{c_{k,n}^2}\left(\mathbb{E}\left[\exp\left(u_n\frac{a_n}{\sum_{k=1}^nq_k}\frac{q_k}{h_k}\frac{\eta_k\left(x\right)}{f\left(x\right)}\right)-1\right]\right)^2\right|}\nonumber\\
&& \leq  \frac{v_n^2}{a_n}\sum_{k=n_0-1}^n
\left(\mathbb{E}\left[\exp\left(u_n\frac{a_n}{\sum_{k=1}^nq_k}\frac{q_k}{h_k}\frac{\eta_k\left(x\right)}{f\left(x\right)}\right)-1\right]\right)^2\nonumber\\
&& = \frac{v_n^2}{a_n}\sum_{k=n_0-1}^n
\left(\mathbb{E}\left[u_n\frac{a_n}{\sum_{k=1}^nq_k}\frac{q_k}{h_k}\frac{\eta_k\left(x\right)}{f\left(x\right)}\right]\right)^2\left(1+o\left(1\right)\right)\nonumber\\
&& = \frac{u^2}{\left(f\left(x\right)\right)^2}a_n\sum_{k=n_0-1}^n
\left(\frac{q_kh_k^{-1}}{\sum_{k=1}^nq_k}\mathbb{E}\left[\eta_k\left(x\right)\right]\right)^2\left(1+o\left(1\right)\right)\nonumber\\
&&=O\left(a_n\frac{\sum_{k=n_0-1}^nq_k^2h_k^4}{\left(\sum_{k=1}^nq_k\right)^2}\right)\nonumber\\
&&=O\left(h_n^5\frac{\sum_{k=n_0-1}^nq_k^2h_k^4}{nq_n^2h_n^4}\left(\frac{nq_n}{\sum_{k=1}^nq_k}\right)^2\right)\nonumber\\
&&=O\left(h_n^5\right)
\end{eqnarray}
which goes to $0$ as $n\to \infty$. Which proves that $\lim_{n\to \infty}\left|R_{n,x}^{\left(2\right)}\left(u\right)\right|=0$. Then, we obtain from (\ref{eq:Moderate}) and (\ref{limbnsum}), $\lim_{n\to \infty}\Lambda_{n,x}\left(u\right)=\Lambda_x^{M}\left(u\right)$. 
$ $\\
{\bf Second case: $\left(v_n\right)\equiv 1$.}
$ $\\
It follows from~(\ref{eq:lamb}) that
\begin{eqnarray}\label{eq:lamb1}
\Lambda_{n,x}\left(u\right)&=&\frac{1}{a_n}\sum_{k=n_0-1}^n\mathbb{E}\left[\exp\left(u\frac{a_n}{\sum_{k=1}^nq_k}\frac{q_k}{h_k}\frac{\eta_k\left(x\right)}{f\left(x\right)}\right)-1\right]\nonumber\\
&&-\frac{1}{2a_n}\sum_{k=n_0-1}^n
\frac{1}{c_{k,n}^2}\left(\mathbb{E}\left[\exp\left(u\frac{a_n}{\sum_{k=1}^nq_k}\frac{q_k}{h_k}\frac{\eta_k\left(x\right)}{f\left(x\right)}\right)-1\right]\right)^2\nonumber\\
&&-\frac{u}{f\left(x\right)} \frac{1}{\sum_{k=1}^nq_k}\sum_{k=n_0-1}^n\frac{q_k}{h_k}\mathbb{E}\left[\eta_k\left(x\right)\right]\nonumber\\
&=&\frac{1}{a_n}\sum_{k=n_0-1}^nh_k\int_{\mathbb{R}^2}\left[\exp\left(\frac{u}{f\left(x\right)}\frac{a_n}{\sum_{k=1}^nq_k}\frac{q_k}{h_k}\left(y-r\left(x\right)\right)K\left(z\right)\right)-1\right]\nonumber\\
&&\times g\left(x,y\right)dzdy\nonumber\\
&&-R_{n,x}^{\left(3\right)}\left(u\right)-R_{n,x}^{\left(4\right)}\left(u\right)+R_{n,x}^{\left(5\right)}\left(u\right)
\end{eqnarray}
with 
\begin{eqnarray*}
R_{n,x}^{\left(3\right)}\left(u\right)&=&\frac{1}{2a_n}\sum_{k=n_0-1}^n
\frac{1}{c_{k,n}^2}\left(\mathbb{E}\left[\exp\left(u\frac{a_n}{\sum_{k=1}^nq_k}\frac{q_k}{h_k}\frac{\eta_k\left(x\right)}{f\left(x\right)}\right)-1\right]\right)^2,\nonumber\\
R_{n,x}^{\left(4\right)}\left(u\right)&=&\frac{u}{f\left(x\right)} \frac{1}{\sum_{k=1}^nq_k}\sum_{k=n_0-1}^n\frac{q_k}{h_k}\mathbb{E}\left[\eta_k\left(x\right)\right],\nonumber\\
R_{n,x}^{\left(5\right)}\left(u\right)&=&\frac{1}{a_n}\sum_{k=n_0-1}^nh_k\int_{\mathbb{R}^2}\left[\exp\left(\frac{u}{f\left(x\right)}\frac{a_n}{\sum_{k=1}^nq_k}\frac{q_k}{h_k}\left(y-r\left(x\right)\right)K\left(z\right)\right)-1\right]
\\
&&\times \left[g\left(x-zh_k,y\right)-g\left(x,y\right)\right]dzdy.
\end{eqnarray*}
It follows from~(\ref{eq:R2term2}), that $\lim_{n\to \infty}\left|R^{\left(3\right)}_{n,x}\left(u\right)\right|=0$. \\
Moreover, in view of~(\ref{limbnsum}) and~(\ref{eq:Eeta}), we have
\begin{eqnarray*}
\left|R_{n,x}^{\left(4\right)}\left(u\right)\right|&=&O\left(\frac{1}{\sum_{k=1}^nq_k}\sum_{k=1}^nq_kh_k^2\right)\\
&=&O\left(\frac{nq_n}{\sum_{k=1}^nq_k}\frac{\sum_{k=1}^nq_kh_k^2}{nq_nh_n^2}h_n^2\right)\\
&=&O\left(h_n^2\right)
\end{eqnarray*} 
which goes to $0$ as $n\to \infty$.\\
Let us now prove that
\begin{eqnarray*}
\lim_{n\to\infty}\left|R_{n,x}^{\left(5\right)}\left(u\right)\right|=0.
\end{eqnarray*}
Set $M>0$ and $\varepsilon>0$; we then have
\begin{eqnarray*}
R_{n,x}^{\left(5\right)}\left(u\right)
&=&\frac{1}{a_n}\sum_{k=n_0-1}^nh_k\int_{\left\{\left|z\right|\leq M\right\}\times\mathbb{R}}\left[\exp\left(\frac{u}{f\left(x\right)}\frac{a_n}{\sum_{k=1}^nq_k}\frac{q_k}{h_k}\left(y-r\left(x\right)\right)K\left(z\right)\right)-1\right]
\\
&&\times\left[g\left(x-zh_k,y\right)-g\left(x,y\right)\right]dzdy\\
&&+\frac{1}{a_n}\sum_{k=n_0-1}^nh_k\int_{\left\{\left|z\right|> M\right\}\times\mathbb{R}}\left[\exp\left(\frac{u}{f\left(x\right)}\frac{a_n}{\sum_{k=1}^nq_k}\frac{q_k}{h_k}\left(y-r\left(x\right)\right)K\left(z\right)\right)-1\right]
\\
&&\times\left[g\left(x-zh_k,y\right)-g\left(x,y\right)\right]dzdy\\
&=&I+II.
\end{eqnarray*}
Using~(\ref{cte:c}), and since for any $t\in \mathbb{R}$, $\left|e^{t}-1\right|\leq \left|t\right|e^{\left|t\right|}$, we have
\begin{eqnarray*}
\left|II\right|
&\leq& \left|u\right|\sum_{k=n_0-1}^n\frac{q_k}{\sum_{k=1}^nq_k}\int_{\left\{\left|z\right|> M\right\}\times\mathbb{R}}\left|y-r\left(x\right)\right|\left|K\left(z\right)\right|\\
&&\times \exp\left(c\frac{u}{f\left(x\right)}\left|y-r\left(x\right)\right|\left|K\left(z\right)\right|\right)
\left|g\left(x-zh_k,y\right)-g\left(x,y\right)\right|dzdy\\
&\leq&\left|u\right|\sum_{k=n_0-1}^n\frac{q_k}{\sum_{k=1}^nq_k}\int_{\left\{\left|z\right|> M\right\}}\left|K\left(z\right)\right|\\
&&\times\left[\int_{\mathbb{R}}\left|y-r\left(x\right)\right|\exp\left(c\frac{u}{f\left(x\right)}\left|y-r\left(x\right)\right|\left|K\left(z\right)\right|\right)
g\left(x-zh_k,y\right)dy\right]dz\\
&&+\left|u\right|\sum_{k=n_0-1}^n\frac{q_k}{\sum_{k=1}^nq_k}\int_{\left\{\left|z\right|> M\right\}}\left|K\left(z\right)\right|\\
&&\times\left[\int_{\mathbb{R}}\left|y-r\left(x\right)\right|\exp\left(c\frac{u}{f\left(x\right)}\left|y-r\left(x\right)\right|\left|K\left(z\right)\right|\right)
g\left(x,y\right)dy\right]dz\\
&\leq & A \int_{\left\{\left|z\right|> M\right\}}\left|K\left(z\right)\right|dz,
\end{eqnarray*}
where $A$ is a constant; this last inequality follows from~(\ref{eq:bounded}) and from the fact that $K$ is bounded.\\
Now, since $K$ is integrable, we can choose $M$ such that 
\begin{eqnarray*}
\left|II\right|\leq \frac{\varepsilon}{2}.
\end{eqnarray*}
Now, for $I$, we write
\begin{eqnarray*}
I&=&\frac{1}{a_n}\sum_{k=n_0-1}^nh_k\int_{\left\{\left|z\right|\leq M\right\}\times\mathbb{R}}\exp\left(\frac{u}{f\left(x\right)}\frac{a_n}{\sum_{k=1}^nq_k}\frac{q_k}{h_k}\left(y-r\left(x\right)\right)K\left(z\right)\right)
\\
&&\times \left[g\left(x-zh_k,y\right)-g\left(x,y\right)\right]dzdy\\
&&-\frac{1}{a_n}\sum_{k=n_0-1}^nh_k\int_{\left\{\left|z\right|\leq M\right\}\times\mathbb{R}}\left[g\left(x-zh_k,y\right)-g\left(x,y\right)\right]dzdy
\end{eqnarray*}
In view of $\left(M4\right)$,~(\ref{eq:bounded}),~(\ref{limbnsum}), the dominated convergence theorem ensure that both integrals converge to $0$. We deduce that for $n$ large enough,
\begin{eqnarray*}
\left|I\right|\leq \frac{\varepsilon}{2},
\end{eqnarray*}
which ensures that $\lim_{n\to \infty}\left|R_{n,x}^{\left(5\right)}\left(u\right)\right|=0$.\\
Then, it follows from~(\ref{eq:lamb1}), and~(\ref{limbnsum}) and from some analysis considerations that
\begin{eqnarray*}
\lefteqn{\lim_{n\to \infty}\Lambda_{n,x}\left(u\right)}\\
&=&\lim_{n\to \infty}\frac{1}{n}\sum_{k=n_0-1}^n\left(\frac{k}{n}\right)^{-a} \\
&&\times\int_{\mathbb{R}^2}\left[\exp\left(\left(1-q\right)\left(\frac{k}{n}\right)^{a-q}\frac{u}{f\left(x\right)}\left(y-r\left(x\right)\right)K\left(z\right)\right)-1\right]g\left(x,y\right)dzdy\nonumber\\
&=&\left(1-q\right)\int_{\left[0,1\right] \times \mathbb{R}^2}s^{-a}\left(\exp\left(us^{a-q}K\left(z\right)\frac{\left(y-r\left(x\right)\right)}{f\left(x\right)}\right)-1\right)g\left(x,y\right)dsdzdy\\
&=&\Lambda_{x}^{L}\left(u\right)
\end{eqnarray*}
and thus Lemma 1 is proved.
\end{proof}
\subsection{Proof of Proposition~\ref{prop:LMDP}}\label{proof:prop2}
\begin{proof}
$ $\\
To prove Proposition~\ref{prop:LMDP}, we apply Proposition~\ref{prop:convexe}, Lemma~\ref{lemma:convLam} and the following result (see A.A. Puhalskii, \cite{Puh94}).
\begin{lemma}\label{lemma:puhalskii}
Let $\left(Z_n\right)$ be a sequence of real random variables, $\left(\nu_n\right)$ a positive sequence satisfying $\lim_{n\to \infty}\nu_n=+\infty$, and suppose that there exists some convex non-negative function $\Gamma$ defined on $\mathbb{R}$ such that
\begin{eqnarray*}
\forall u\in \mathbb{R}, \lim_{n\to \infty}\frac{1}{\nu_n}\log \mathbb{E}\left[\exp\left(u\nu_nZ_n\right)\right]=\Gamma\left(u\right).
\end{eqnarray*}
If the Legendre function $\Gamma^*$ of $\Gamma$ is a strictly convex function, then the sequence $\left(Z_n\right)$ satisfies a \texttt{LDP} of speed $\left(\nu_n\right)$ and good rate fonction $\Gamma^*$.
\end{lemma}
In our framework, when $v_n\equiv 1$, we take $Z_n=\overline \rho_n\left(x\right)-\mathbb{E}\left(\overline \rho_n\left(x\right)\right)$, $\nu_n=nh_n$ with $h_n=cn^{-a}$ where $c>0$ and $a\in\left]1-\alpha,\left(4\alpha-3\right)/2\right[$ (with $\alpha\in ]\frac{3}{4},1]$), and the weight $\left(q_n\right)=\left(c^{\prime}n^{-q}\right)$ with $c^{\prime}>0$ and $q<\min\left\{1-2a,\left(1+a\right)/2\right\}$, and $\Gamma=\Lambda_x^{L}$. In this case, the Legendre transform of $\Gamma=\Lambda_x^{L}$ is the rate function $I_{a,q,x}\left(t\right)$ which is strictly convex by Proposition~\ref{prop:convexe}. Otherwise, when, $v_n\to \infty$, we take $Z_n=v_n\left(\overline \rho_n\left(x\right)-\mathbb{E}\left[\overline \rho_n\left(x\right)\right]\right)$, $\nu_n=nh_n/v_n^2$ and $\Gamma=\Lambda_x^M$; $\Gamma^*$ is then the quadratic rate function $J_{a,q,x}$ defined in~(\ref{MDPalgoreg}) and thus Proposition~\ref{prop:LMDP} follows.
\end{proof}
\subsection{Proof of Proposition~\ref{prop:Bn}}\label{proof:prop3}
\begin{proof}
$ $\\
It follows from~(\ref{rhon}),~(\ref{av:Rn0}),~(\ref{limbnsum}) and~(\ref{eq:Eeta}), that
\begin{eqnarray*}
\mathbb{E}\left[\overline \rho_n\left(x\right)\right]-r\left(x\right)&=&\frac{1}{f\left(x\right)}\mathbb{E}\left[T_n\left(x\right)\right]\\
&=&\frac{1}{f\left(x\right)}\frac{1}{\sum_{k=1}^nq_k}\sum_{k=n_0-1}^n\frac{q_k}{h_k}\mathbb{E}\left[\eta_k\left(x\right)\right]\\
&=&\frac{1}{f\left(x\right)}\frac{\sum_{k=n_0-1}^nq_kh_k^2}{\sum_{k=1}^nq_k}m^{\left(2\right)}\left(x\right)f\left(x\right)\left[1+o\left(1\right)\right]\\
&=&h_n^2\frac{1-q}{1-q-2a}m^{\left(2\right)}\left(x\right)\left[1+o\left(1\right)\right]\\
&=&O\left(h_n^2\right).
\end{eqnarray*}
\end{proof}
\subsection{Proof of Proposition~\ref{prop:convexe}}\label{proof:prop1}
\begin{proof}
$ $\\
\begin{itemize}
\item Since $\left|e^t-1\right|\leq \left|t\right|e^{\left|t\right|}$ $\forall t\in \mathbb{R}$, it follows from~(\ref{eq:bounded}) and $\left(L1\right)$, that
\begin{eqnarray*}
\left|\psi_{a,q,x}\left(u\right)\right|&\leq& \left(1-q\right)\int_{\left[0,1\right]\times\mathbb{R}^2}s^{-a}\left|\exp\left(us^{a-q}K\left(z\right)\frac{\left(y-r\left(x\right)\right)}{f\left(x\right)}\right)-1\right|g\left(x,y\right)dsdzdy\\
&\leq& \left(1-q\right)\frac{\left|u\right|}{f\left(x\right)}\int_{\left[0,1\right]\times\mathbb{R}^2}s^{-q}\left|y-r\left(x\right)\right|\left|K\left(z\right)\right|\\
&&\times \exp\left(\left|u\right|\frac{\left|y-r\left(x\right)\right|}{f\left(x\right)}\|K\|_{\infty}\right)g\left(x,y\right)dsdzdy\\
&\leq& \frac{\left|u\right|}{f\left(x\right)}\int_{\mathbb{R}^2}\left|y-r\left(x\right)\right|\left|K\left(z\right)\right|\exp\left(\left|u\right|\frac{\left|y-r\left(x\right)\right|}{f\left(x\right)}\|K\|_{\infty}\right)g\left(x,y\right)dzdy\\
&=& \frac{\left|u\right|}{f\left(x\right)}\int_{\mathbb{R}}\left|K\left(z\right)\right|dz\int_{\mathbb{R}}\left|y-r\left(x\right)\right|\exp\left(\left|u\right|\frac{\left|y-r\left(x\right)\right|}{f\left(x\right)}\|K\|_{\infty}\right)g\left(x,y\right)dy\\
&<&\infty
\end{eqnarray*}
which ensures the existence of $\psi_{a,q,x}$. It is straightforward to check that $\psi_{a,q,x}$ is twice differentiable, with
\begin{eqnarray*}
\psi_{a,q,x}^{\prime}\left(u\right)&=&\left(1-q\right)\int_{\left[0,1\right]\times\mathbb{R}^2}s^{-q}K\left(z\right)\frac{\left(y-r\left(x\right)\right)}{f\left(x\right)}\\
&&\times 
\exp\left(us^{a-q}K\left(z\right)\frac{\left(y-r\left(x\right)\right)}{f\left(x\right)}\right)g\left(x,y\right)dsdzdy\\
\psi_{a,q,x}^{\prime\prime}\left(u\right)&=&\left(1-q\right)\int_{\left[0,1\right]\times\mathbb{R}^2}s^{a-2q}\left(K\left(z\right)\right)^2\left(\frac{\left(y-r\left(x\right)\right)}{f\left(x\right)}\right)^2\\
&&\times
\exp\left(us^{a-q}K\left(z\right)\frac{\left(y-r\left(x\right)\right)}{f\left(x\right)}\right)g\left(x,y\right)dsdzdy.
\end{eqnarray*}
Since $\psi_{a,q,x}^{\prime\prime}\left(u\right)>0$ $\forall u\in \mathbb{R}$, $\psi^{\prime}_{a,q,x}$ is increasing on $\mathbb{R}$, and $\psi_{a,q,x}$ is strictly convex on $\mathbb{R}$. It follows that its Cramer transform $I_{a,q,x}$ is a good rate function on $\mathbb{R}$ (see A. Dembo and O. Zeitouni \cite{Dem98}) and $\left(i\right)$ of Proposition~\ref{prop:convexe} is proved.
\item Let us now assume that $\lambda\left(O_-\right)=0$. We then have
\begin{eqnarray*}
\lim_{u\to-\infty}\psi_{a,q,x}^{\prime}\left(u\right)=0\quad \mbox{and}\quad \lim_{u\to+\infty}\psi_{a,q,x}^{\prime}\left(u\right)=+\infty
\end{eqnarray*}
so that the range of $\psi_{a,q,x}^{\prime}$ is $]0,+\infty[$. Moreover 
\begin{eqnarray*}
\lim_{u\to-\infty}\psi_{a,q,x}\left(u\right)&=&\left\{\begin{array}{lllll}
-\left(1-q\right)/\left(1-a\right)\lambda\left(S_+\right)f\left(x\right) & \mbox{if} & \lambda\left(S_+\cap T_+\right)>0\\
-\left(1-q\right)/\left(1-a\right)\lambda\left(S_-\right)f\left(x\right) & \mbox{if} & \lambda\left(S_-\cap T_-\right)>0\\
\end{array}
\right.\\
\end{eqnarray*}
(which can be $-\infty$). This implies in particular that 
\begin{eqnarray*}
I_{a,q,x}\left(0\right)&=&\left\{\begin{array}{lllll}
\left(1-q\right)/\left(1-a\right)\lambda\left(S_+\right)f\left(x\right) & \mbox{if} & \lambda\left(S_+\cap T_+\right)>0\\
\left(1-q\right)/\left(1-a\right)\lambda\left(S_-\right)f\left(x\right) & \mbox{if} & \lambda\left(S_-\cap T_-\right)>0\\
\end{array}
\right.\\
\end{eqnarray*}
Now, when $t < 0$, $\lim_{u\to-\infty}\left(ut-\psi_{a,q,x}\left(u\right)\right)=+\infty$ and $I_{a,q,x}\left(t\right)=+\infty$. Since $\psi^{\prime}_{a,q,x}$ is increasing with range $]0,+\infty[$, when $t>0$, $\sup_u\left(ut-\psi_{a,q,x}\left(u\right)\right)$ is reached for $u_0\left(t\right)$ such that $\psi_{a,q,x}\left(u_0\left(t\right)\right)=t$, i.e. for $u_0\left(t\right)=\left(\psi_{a,q,x}^{\prime}\right)^{-1}\left(t\right)$; this prove~(\ref{eq:taurev}). (Note that, since $\psi_{a,q,x}^{\prime \prime}\left(t\right)>0$, the function $t\mapsto u_0\left(t\right)$ is differentiable on $]0,+\infty[$). Now, differentiating~(\ref{eq:taurev}), we have
\begin{eqnarray*}
I^{\prime}_{a,q,x}\left(t\right)&=&u_0\left(t\right)+tu_0^{\prime}\left(t\right)-\psi_{a,q,x}^{\prime}\left(u_0\left(t\right)\right)u_0^{\prime}\left(t\right)\\
&=&\left(\psi_{a,q,x}^{\prime}\right)^{-1}\left(t\right)+tu_0^{\prime}\left(t\right)-tu_0^{\prime}\left(t\right)\\
&=&\left(\psi_{a,q,x}^{\prime}\right)^{-1}\left(t\right).
\end{eqnarray*}
Since $\left(\psi_{a,q,x}^{\prime}\right)^{-1}$ is an increasing function on $]0,+\infty[$, it follows that $I_{a,q,x}$ is strictly convex on $]0,+\infty[$ (and differentiable). Thus (ii) is proved.\\
\item 
We Assume that $\lambda\left(O_-\right)>0$. In this case, $\psi_{a,q,x}^{\prime}$ can be rewritten as
\begin{eqnarray*}
\psi_{a,q,x}^{\prime}\left(u\right)&=&\left(1-q\right)\int_{\left[0,1\right]\times\left(\mathbb{R}^2\cap O_+\right)}s^{-q}K\left(z\right)\frac{\left(y-r\left(x\right)\right)}{f\left(x\right)}\\
&&\times
\exp\left(us^{a-q}K\left(z\right)\frac{\left(y-r\left(x\right)\right)}{f\left(x\right)}\right)g\left(x,y\right)dsdzdy\\
&&+\left(1-q\right)\int_{\left[0,1\right]\times\left(\mathbb{R}^2\cap O_-\right)}s^{-q}K\left(z\right)\frac{\left(y-r\left(x\right)\right)}{f\left(x\right)}\\
&&\times 
\exp\left(us^{a-q}K\left(z\right)\frac{\left(y-r\left(x\right)\right)}{f\left(x\right)}\right)g\left(x,y\right)dsdzdy
\end{eqnarray*}
and we have
\begin{eqnarray*}
\lim_{u\to -\infty} \psi^{\prime}_{a,q,x}\left(u\right)=-\infty\quad \mbox{and}\quad \lim_{u\to +\infty} \psi^{\prime}_{a,q,x}\left(u\right)=+\infty
\end{eqnarray*}
so that the range of $\psi^{\prime}_{a,q,x}$ is $\mathbb{R}$ in this case. The proof of $\left(iii\right)$ follows the same lines as previously, except that, in the present case, $\left(\psi^{\prime}_{a,q,x}\right)^{-1}$ is defined on $\mathbb{R}$, and not only on $]0,+\infty[$.
\end{itemize}
\end{proof}


\begin{thebibliography}{99}
\bibitem{Bo73}
\uppercase{Bojanic, R.---Seneta, E.}:
\textit{A unified theory of regularly varying sequences},
Math. Z. {\bf 134} (1973),~91--106.

\bibitem{Dem98} 
\uppercase{Dembo, A.---Zeitouni, O.}:
\textit{Large deviations techniques and applications}, Springer,
Applications of mathematics, New-York, 1998.


\bibitem{Ga73}
\uppercase{Galambos, J.---Seneta, E.}:
\textit{Regularly varying sequences}, Proc.
Amer. Math. Soc. {\bf 41} (1973),~110--116.

\bibitem{Jot06} 
\uppercase{Joutard, C.}: 
\textit{Sharp large deviations in nonparametric estimation}, J. of Nonparam. Stat. {\bf 18} (2006),~293--306.

\bibitem{Lou99} 
\uppercase{Louani, D.}:
\textit{Some large deviations limit theorems in conditionnal nonparametric statistics}, Statistics. {\bf 33} (1999),~171--196.

\bibitem{Nad64} 
\uppercase{Nadaraya, E. A.}:
\textit{On estimating regression}, Theory Probab. Appl. {\bf 10} (1964),~186--190.

\bibitem{Mok07} 
\uppercase{Mokkadem, A.---Pelletier, M.}:
\textit{A companion for the Kiefer-Wolfowitz-Blum stochastic approximation algorithm}, Ann. Statist. {\bf 35} (2007),~1749--1772.

\bibitem{Mok08} 
\uppercase{Mokkadem, A.---Pelletier, M.---Thiam, B.}:
\textit{Large and moderate deviations principles for kernel estimators of the multivariate regression}, Mathematical Methods of Statistics. {\bf 17} (2008),~1--27.


\bibitem{Mok09} 
\uppercase{Mokkadem, A.---Pelletier, M.---Slaoui, Y.}:
\textit{Revisiting R\'ev\'esz's stochastic approximation method for the estimation of a regression function}, ALEA Lat. Am. J. Probab. Math. Stat. {\bf 6} (2009),~63--114.

\bibitem{Puh94} 
\uppercase{Puhalskii, A. A.}:
\textit{The method of stochastic exponentials for large deviations}, Stochastic Process. Appl. {\bf 54} (1994),~45--70.

\bibitem{Rev73}
\uppercase{R\'ev\'esz, P.}:
\textit{Robbins-Monro procedure in a Hilbert space and its
application in the theory of learning processes I}, Studia Sci. Math. Hung. {\bf 8} (1973),~391--398.

\bibitem{Rev77}
\uppercase{R\'ev\'esz, P.}:
\textit{How to apply the method of stochastic approximation in the non-parametric estimation of a regression function}, Math. Operationsforsch. Statist., Ser. Statistics. {\bf 8} (1977),~119--126.

\bibitem{Tsy90} 
\uppercase{Tsybakov, A.B.}:
\textit{Recurrent estimation of the mode of a multidimensional distribution}, Problems Inform. Transmission. {\bf 26} (1990)~31--37.

\bibitem{Wat64} 
\uppercase{Watson, G. S.}:
\textit{Smooth regression analysis}, Sankhya Ser. A. {\bf 26} (1990)~359--372.

\end{thebibliography}
\end{document}